\documentclass[conference]{IEEEtran}

\IEEEoverridecommandlockouts 

\usepackage{cite}

\ifCLASSINFOpdf
   \usepackage[pdftex]{graphicx}
   \graphicspath{{./Images/}}

   \DeclareGraphicsExtensions{.pdf,.jpeg,.png}
\else

\fi

\usepackage[cmex10]{amsmath}

\usepackage{amssymb}
\usepackage{amsthm}
\usepackage{amscd}

\usepackage{tikz}
\usetikzlibrary{shapes,arrows}
\usepackage{multirow}
\usepackage{booktabs,colortbl}
\begin{document}
%
\title{A robust convex optimization framework for autonomous network planning under load uncertainty}

\author{%

\IEEEauthorblockN{Beno\^it MARTIN*, Emmanuel DE JAEGER}
\IEEEauthorblockA{Mechatronic, Electrical Energy, and Dynamic Systems\\ (MEED)\\
 iMMC,UCL\\
Louvain-la-Neuve, Belgium\\
\{benoit.martin,emmanuel.dejaeger\}@uclouvain.be}
\and
\IEEEauthorblockN{Fran\c cois GLINEUR}
\IEEEauthorblockA{Center for Operations Research and Econometrics\\ (CORE)\\
 UCL\\
Louvain-la-Neuve, Belgium\\
{francois.glineur@uclouvain.be}
}

%

}


\maketitle

\begin{abstract}
Autonomous microgrid planning is a Mixed-Integer Non Convex decision problem that requires to consider investments in both distribution and generation capacity and represents significant computation challenges. We proposed in a previous publication a deterministic Second-Order Cone (SOC) relaxation of this problem that made it computationally tractable for real-size cases. However,  this problem is subject to considerable uncertainty emanating from load consumption, RES-based generation and contingencies. In this paper, we thus present a robust optimization approach that extends our previous work by including load related uncertainty at the cost of a substantial increase of the computational burden. The results show that significantly higher investment and operational costs are incurred to account for the load related uncertainty\\
\end{abstract}


\section{Introduction}
Electrification is considered as a major developing factor in modern societies. However, even if 84\% of world population did have access to electricity in 2016, this figure hides significant disparities. First, between countries, as this percentage dropped to 19\% for sub-Saharan Africa during the same year. Then, between areas as the vast majority of people without electricity access (around 80\%) lived in rural zones worldwide \cite{WEO_IEA}.\\
Rural electrification in developing countries is a significant challenge. As a matter of fact, traditional extension of centralized grid may be inefficient in this context for reasons such as capital scarcity, remoteness and lack of reliability \cite{schnitzer_microgrids_????}.
Autonomous microgrids thus offer an efficient alternative for rural electrification as they are less capital-intensive and offer better reliability thanks to distributed energy resources.\\
Autonomous microgrid planning consists in making investment decisions concerning an isolated microgrid on a predefined planning horizon. The isolated (or autonomous) character of such power systems implies that sufficient generation capacity should be placed to meet the total demand occurring in the system. Considering a set of $n$ nodes representing future consumption points to electrify and consumption profiles for these nodes, the problem consists in answering the following questions in such a way as to minimize the total cost (OPEX and CAPEX) of the system on the planning horizon:
\begin{itemize}
\item Which distribution and generation assets should be placed?
\item Where to place them?
\item When to place them?
\end{itemize}
The problem is inherently uncertain as consumption profiles are forecasts subject to errors. Furthermore, if RES based generators are considered, their power output is also uncertain. At last, uncertainty also arises from microgrid components (lines, generators) subject to contingencies. \\In this paper, we present a robust approach to autonomous microgrid planning considering load consumption uncertainty. This paper extends our previous work on a deterministic formulation for this problem \cite{martin_comparison_????} with the computation of worst operating scenarios developed in \cite{capitanescu_computation_2013}.\\ 
This paper is organized as follows: the deterministic formulation previously developed is presented in section \ref{DetFor},  section \ref{UncMod} describes the uncertainty model that is used and section \ref{ProbScen} presents the approach used to compute problematic operating scenarios. We present the results in section \ref{Res}  before concluding in section \ref{Concl}.
\section{Deterministic formulation of autonomous microgrid planning}
\label{DetFor}

Autonomous microgrid planning is a high-dimensional problem with a lot of discrete investment decision variables. Firstly, generation and distribution capacity have to be installed, we thus consider them simultaneously in the planning problem as sequential planning of these two elements would be suboptimal. Secondly, autonomous microgrid planning differs from traditional expansion planning in the way that microgrids are often built from scratch in this context, which necessitates a lot of decision variables. The joint consideration of large number of discrete decision variables is known to lead to combinatorial explosion.  Another salient feature of the problem is its non-convexity which is caused by power flow equations.\\
We presented in a previous work \cite{martin_comparison_????} four different convex formulations of the planning problem that could be cast as mixed-integer convex programs. The convexity of these formulations allows to reach a single global optimum. The first three formulations are linear approximations of the non-convex problem while the last one is a second-order cone relaxation of the original problem. In the latter case, the objective value corresponding to the relaxation global optimum is thus a lower bound of the optimal objective value for the non-convex problem. For this work, Benders decomposition has been successfully used to manage the dimensionality of the integer part of the problem.

Our goal is to take investment decisions (lines, generators) and operational decisions (production of generators) in order to minimize the total net present value of the system on the planning horizon . Investment decisions are taken once a year while operational decisions are taken once an hour. We simulate a limited amount of representative days per year.

\noindent We consider a set of \textit{n} nodes representing consumption points to electrify. Data for these nodes includes location (coordinates) and hourly consumption profiles for two typical days of the year. Consumptions are assumed to increase at a uniform rate each year, triggering the need for reinforcements. Available investment options consist of diesel generator sets with a linear operational cost function and overhead cables. We consider a unique size available for the lines and the generators. However, while there may be at most one generator installed at a node \textit{i}, there may be several lines placed in parallel between nodes \textit{i }and \textit{j }which is equivalent to a bigger line. 

\noindent The objective is to minimize the NPV of the system, which is the sum of discounted yearly cash flows (CAPEX and OPEX).  In these expressions, ${\omega }_{ijy}\ $is a binary variable equal to 1 if there is at least one conductor between nodes \textit{i} and \textit{j }at year \textit{y}, ${\gamma }_{ijy}\ is$ the number of lines in parallel between nodes \textit{i }and \textit{j }at year \textit{y,}$\ {\sigma }_{iy}\ $is a binary variable equal to 1 if there is a generator of fixed size at node \textit{i }at year \textit{y,}$\ P_{Git}$\textit{ }is the active power produced at node \textit{i }at period \textit{t.} The parameters are the following:\textit{ }$C_{cond}$\textit{  }is the cost of a single conductor (\$/km), $C_{pole}$ is the cost of poles (\$/km) (a unique pole is required regardless of the number of conductors),$\ D_{ij}$ is the distance between nodes \textit{i} and \textit{j }(km),\textit{ }  $C_{Gen}$ is the fixed cost for installing a generator, \textit{a }and \textit{b} are the parameters of the generator linear cost function and \textit{ra} is the discount rate. As we only simulate a limited amount of days per year, we multiply the fuel costs for these days by a suitable scaling factor $H$ to represent the yearly operation cost.

\begin{flalign}
&CAPEX_{Dist,y}=\sum_{\left(i,j\right)}{\left({\gamma }_{ijy}-{\gamma }_{ijy-1}\right){D_{ij}C}_{cond}}&
\end{flalign}

\begin{flalign*}
+\left({\omega }_{ijy}-{\omega }_{ijy-1}\right){D_{ij}C}_{pole}&
\end{flalign*}
\begin{flalign}
&CAPEX_{Gen,y}=\sum_i{\left({\sigma}_{iy}-{\sigma}_{iy-1}\right)C_{Gen}}&
\end{flalign}
\begin{flalign}
&OPEX_y=H\sum_i{\sum_{t\in y}{(a{\sigma}_{iy}+bP_{Git}})}&
\end{flalign}
\begin{flalign}
&NPV=\sum^Y_{y=1}{\frac{1}{{\left(1+ra\right)}^y}[CAPEX_{Dist,y}+CAPEX_{Gen,y}+OPEX_y]}&
\end{flalign}

Eq. (5) express the fact that ${\gamma }_{ijy}$ shouldn't exceed the maximal amount of parallel lines $\overline{\xi }$ and that investments are permanent (i.e. cannot be unmade in following years). Eq. (6) also states that investments in generation can't be unmade as ${\sigma }_{iy}$, should be increasing though time. Eq. (7) simply expresses the symmetry of the problem regarding lines. 

\begin{flalign}
\gamma_{ijy-1}\le {\gamma }_{ijy}\le \overline{\xi}&
\end{flalign}
\begin{flalign}  
{\sigma }_{iy-1}\le {\sigma }_{iy} &
\end{flalign}
\begin{flalign}
\omega_{ijy}=\omega_{jiy}, \gamma_{ijy}=\gamma_{jiy}&
\end{flalign}

Eqs. (8) to (11) force the network to be at least radial (and allows it to be meshed) where \textit{n} is the number of nodes and $f_{jiy}$ is a fictitious flow only used to ensure connectivity of the network. The idea behind those constraints is to have a fictitious source supplying \textit{(n-1)} units at node 1 and fictitious sinks at other nodes that consume \textit{1}. As there is only one source, eqs. (8) to (11) ensure that there is no island in the network.

\begin{flalign} 
\sum_{\left(i,j\right)\ }{{\omega }_{ijy}\ge 2\times \left(n-1\right)}&
\end{flalign}
\begin{flalign} 
f_{ijy}\le {\omega }_{ijy}\times n &
\end{flalign}
\begin{flalign} 
\sum_{(i,j)}{f_{1jy}=n-1} &
\end{flalign}
\begin{flalign}
\sum_{(i,j)}{f_{jiy}=1+\ }\ \sum_{(i,j)}{f_{ijy}} &
\end{flalign}

Eq. (12) ensures that the active power produced at node \textit{i }is smaller than the generation capacity installed at this node and larger than the technical minimum. Eq. (13) represents reactive capabilities of generators, with ${\mathrm{cos} \left(\mathrm{\Phi }\right)\ }$ being the minimal power factor (reactive or inductive) of generation units.

\begin{flalign}
{\sigma \ }_{iy}\ \underline{P}\ \le P_{Git}\le {\sigma }_{iy}\ \overline{P} &
\end{flalign}
\begin{flalign*}
-P_{Git}\times \mathrm{tan}\mathrm{}({{\mathrm{cos}}^{-1} \left({\mathrm{cos} \left(\mathrm{\Phi }\right)\ }\right)\ }\le Q_{Git}&
\end{flalign*}
\begin{flalign}
{\ \ Q}_{Git}\le P_{Git}\times \mathrm{tan}\mathrm{}({{\mathrm{cos}}^{-1} \left({\mathrm{cos} \left(\mathrm{\Phi }\right)\ }\right)\ } & 
\end{flalign}

Eqs. (14) and (15) represent the active and reactive nodal power balance respectively, $P_{Cit}$ and $Q_{Cit}$ representing active and reactive power consumptions at node \textit{i }at period \textit{t }while$\ p_{ijt}$ and $q_{ijt}\ $represent active and reactive power flows from \textit{i }to \textit{j} at period \textit{t}.

\begin{flalign}
P_{Git}-P_{Cit}=\sum_{\left(i,j\right)\ }{p_{ijt}} &
\label{active_balance}
\end{flalign}
\begin{flalign}
Q_{Git}-Q_{Cit}=\sum_{(i,j)\ }{q_{ijt}} &
\label{reactive_balance}
\end{flalign}

\noindent We define binary variables $loi_{ijyk}$ to be equal to 1 if the amount of parallel lines between \textit{i }and \textit{j }is greater or equal to \textit{k }and zero otherwise, which is expressed by eqs. (16) and (17). These variables are used to write constraints (18), (19), (22) and (23) for each possible level of investment in lines such that we avoid bilinear terms.

\begin{flalign}
\sum^{\overline{\xi}}_{k=1}{loi_{ijky}={\gamma}_{ijy}}&
\end{flalign}
\begin{flalign}
{\omega}_{ijy}=loi_{ijy1}& 
\end{flalign}

The last constraints express the physics of power flows. $\mathrm{\Psi }_{ijt}$ represents the squared amplitude of line current and $\nu_{it}$ represents the squared voltage amplitude. Eqs. (18) expresses active losses and is written such that the only active constraint is the one corresponding to the actual amount of parallel lines between \textit{i}and \textit{j} (i.e. to the unique k such that $loi_{ij,k,y}=1$ and  $loi_{ij,k+1,y}=0$),  in order to avoid bilinear terms. The reactive losses on the line are similarly developed by replacing $p_{ijt}$ and $r$ by $q_{ijt}$ and $x$ respectively in eq. 18. Parameters $\textit{r}$ and $\textit{x}$ are the line resistance and reactance per unit length and $\textit{M}_1$ is a large enough constant. Finally, eqs. (19) force active and reactive losses to be positive on every line. While being redundant, these constraints considerably tighten the resulting model.

\noindent 

\begin{flalign*}
-\left(1-(loi_{ijky}-loi_{ijk+1y})\right)M_1&
\end{flalign*}
\begin{flalign*}
\le p_{ijt}+p_{jit}-\frac{rD_{ij}}{k}{\mathrm{\Psi }}_{ijt}&
\end{flalign*}
\begin{flalign}
\le \left(1-(loi_{ijky}-loi_{ijk+1y})\right)M_1&
\end{flalign}
\begin{flalign}
p_{ijt}+p_{jit}\ge 0,\ q_{ijt}+q_{jit}\ge 0&
\end{flalign}

Eq. (20)  expresses the fact that power flowing in a line is the product of node voltage and line current. It is relaxed as an inequality and has the form of a (convex) rotated second order cone constraint. 

\noindent 

\begin{flalign}
p^2_{ijt}+q^2_{ijt}\le {\mathrm{\Psi }}_{ijt}{\nu }_{it}&
\end{flalign}

Eq. (21) expresses voltage drops and is written in a way similar to (18).  Eq. (23) expresses nodal voltage bounds.

\begin{flalign*}
-\left(1-(loi_{ijky}-loi_{ijk+1y})\right)M_2&
\end{flalign*}
\begin{flalign*}
\le {\nu }_{jt}-{\nu }_{it}+2D_{ij}(\frac{r}{k}p_{ijt}+\frac{x}{k}q_{ijt}-\frac{D_{ij}}{k^2}\left(r^2+x^2\right){\mathrm{\Psi }}_{ijt})
\end{flalign*}
\begin{flalign}
\le \left(1-(loi_{ijky}-loi_{ijk+1y})\right)M_2&
\end{flalign}
\begin{flalign}
{\underline{v}}^2\le {\nu }_{it}\le {\overline{v}}^2&
\end{flalign}

Eq. (23) is the line thermal rating constraint. It is a SOC constraint.

\begin{flalign}
p^2_{ijt}+q^2_{ijt}\le {\gamma }^2_{ijy}{\overline{S}}^2&
\end{flalign}

Finally, as proposed in [4], we introduce lower and upper bounds on voltage angle differences even if this formulation doesn't include angles explicitly. As a matter of fact, these constraints significantly tighten the model. For brevity, we only present the general form of these constraints using the set of available variables. The bilinear term ${\gamma }_{ijy}{\nu }_{it}\ $can be replaced by an appropriate lift-and-project relaxation and ``big M'' constraints similar to (18) and (21). The parameter ${\theta }^{\mathrm{\Delta }}\ $is the maximum angle difference allowed between two nodes.

\begin{flalign*}
rD_{ij}\left(q_{ijt}+{\mathrm{tan} \left({\theta }^{\mathrm{\Delta }}\ \right)\ }p_{ijt}\right)+xD_{ij}\left({\mathrm{tan} \left({\theta }^{\mathrm{\Delta }}\ \right)\ }q_{ijt}-p_{ijt}\right)&
\end{flalign*}
\begin{flalign}
\le {\mathrm{tan} \left({\theta }^{\mathrm{\Delta }}\ \right)\ }{\nu }_{it}{\gamma }_{ijy}&
\end{flalign}
\begin{flalign*}
xD_{ij}\left(p_{ijt}+{\mathrm{tan} \left({\theta }^{\mathrm{\Delta }}\ \right)\ }q_{ijt}\right)+rD_{ij}\left({\mathrm{tan} \left({\theta }^{\mathrm{\Delta }}\ \right)\ }p_{ijt}-q_{ijt}\right)&
\end{flalign*}
\begin{flalign}
\le {\mathrm{tan} \left({\theta }^{\mathrm{\Delta }}\ \right)\ }{\nu }_{it}{\gamma }_{ijy}&
\end{flalign}

\textbf{}

\noindent This model is computationally intractable in its Mixed Integer Second Order Cone (MISOC) form. We thus apply the Ben-Tal Nemirovski (BTN) relaxation [7] to the SOC constraints (Eqs. (20) and (23)). It consists in replacing the second order cones by cutting places in an efficient way, with arbitrary accuracy. This formulation is named MISOC BTN hereafter.

\noindent 

\noindent 

\noindent 

\noindent 

\noindent

\section{Modeling of uncertainty in autonomous microgrid planning}
\label{UncMod}
Three sources of uncertainty can be distinguished in the autonomous microgrid planning problem: load forecast errors, RES-based generation forecast errors and contingencies. In this paper, we only consider load uncertainty that will be modelled with a rectangular uncertainty set $\Omega =\{\boldsymbol{\omega}\in \mathbb{R}^{n^{\Omega}} \; :\omega_i \in [\omega_i^L;\omega_i^U] \; \forall i \in 1,...,n^{\Omega}\}$. This means that we only consider the interval in which random load consumptions may vary without making any assumption about the distribution of these random variables.\newline
By considering uncertainty in the problem formulation, the aim is to build a microgrid satisfying all constraints defined in the previous section not only for a single scenario , e.g. the most likely one, but also for every possible realization of random variables, i.e. every $\boldsymbol{\omega}\in\Omega$. Autonomous microgrid planning thus becomes a robust optimization (RO) problem.\newline However, such problems are difficult to solve and are generally NP-hard . Indeed, considering continuous random variables potentially leads to an infinite uncertainty space which in turn leads to an infinite amount of constraints to consider in the RO problem \cite{calafiore_scenario_2006}. In \cite{calafiore_scenario_2006}, the authors propose a finite constraint sampling scheme to overcome this problem. They show that the probability of constraint violations rapidly decreases with the amount of samples. They also provide an upper bound on the amount of samples needed to obtain a predefined level of confidence concerning constraint enforcement, which allows to efficiently solve the problem to arbitrary accuracy. In \cite{margellos_road_2014} and \cite{venzke_convex_2017}, the authors propose another method to reduce to a finite size the set of constraints. They show that for a problem with a polytopic uncertainty set $\Omega$ and convex constraints of the form $g(x)\leq 0$ , the body of these constraints will always be maximal on the vertices of $\Omega$. To enforce such constraints for all $\boldsymbol{\omega}\in\Omega$, it is thus sufficient to enforce them on every vertex of $\Omega$. Nonetheless, even in the simple case where $\Omega$ is a rectangular set, the amount of vertices is equal to $2^{n^{\Omega}}$ which rapidly becomes intractable with a growing $n^{\Omega}$.\newline
Consequently, we adopt the approach developed in \cite{capitanescu_computation_2013} which has been used for security planning under uncertainty in transmission networks \cite{capitanescu_cautious_2012}. This approach consists in computing a subset of the vertices of $\Omega$, i.e. a set of scenarios to incorporate in the RO problem, sufficient to guarantee constraint enforcement on the whole uncertainty set $\Omega$. The approach, described in section \ref{ProbScen}, is based on the successive and iterative  computation of an adversarial problem where the infeasibility ( i.e. violation) of the constraints is maximized in order to find problematic scenarios to add to the RO problem and a corrective problem where we try to remove these infeasibilities thanks to remedial actions. In this paper, we consider two sorts of infeasibilities: insufficient generation capacity and line thermal rating violation.

\section{Determination of problematic scenarios}
\label{ProbScen}
The scenario generation algorithm developed in \cite{capitanescu_cautious_2012} can be summarized as follows, $\mathcal{S}$ and $\mathcal{PS}$ being the set of scenarios to consider in the RO problem and the current set of problematic scenarios respectively. All these steps are described in the following subsections.

\begin{enumerate}
\item Initialize $\mathcal{S}$ with the scenario corresponding to the deterministic case
\item Unfix investment variables, solve main problem on $\mathcal{S}$ and then fix investement variables to their current optimal values
\item Reinitialize $\mathcal{PS} \leftarrow \emptyset  $ and solve adversarial problem to compute the current set $\mathcal{PS}$
\item Solve corrective problem $\forall s \in \mathcal{S}$. If there are no more infeasibilities  for $s^*$, then it is not a problematic scenario: $\mathcal{PS} \leftarrow \mathcal{PS}\setminus \{s^*\}$
\item If $\mathcal{PS} = \emptyset$, END. Else, update $\mathcal{S}\leftarrow\mathcal{S} \cup \mathcal{PS}$ and go back to step 2).
\end{enumerate}

We now express active and reactive consumptions as random variables $\tilde{P}_{cit}$ and $\tilde{Q}_{cit}$ in the adversarial problem while they were parameters in the deterministic formulation of section \ref{DetFor}. We thus have two uncertainty sets: $\Omega^{P}$ and $\Omega^{Q}$ for active and reactive power consumptions respectively. $n^{\Omega}$ is equal to $n \times T$ as there is a power consumption forecast for every node and every timestep of the planning horizon. A scenario $s$ thus consists of two matrices $\mathbf{\tilde{P}_{c}}$ and $\mathbf{\tilde{Q}_{c}}$ $\in \mathbb{R}^{n\times T}$ corresponding to a particular realization of random variables.

\subsection{Main problem}

The main problem is the deterministic problem described in section \ref{DetFor} with the following differences: operational variables $P_{Gits},Q_{Gits},p_{ijts},q_{ijts}, \Psi_{ijts}$ and $\nu_{its}$ and operational constraints (12)-(15) and (18)-(25) are now indexed on the scenario set $\mathcal{S}$ as well and the total OPEX is now the expected value of OPEX over all scenarios considering that they all are equiprobable.\newline It is thus a deterministic problem where operational constraints are replicated for each $s \in \mathcal{S}$. This problem is solved at each iteration of the algorithm. It should be emphasized that, while operational variables may be now adapted for each scenario, investment variables remain common to every scenario in order to find a unique investment plan suitable for the whole set $\mathcal{S}$.

\subsection{Adversarial problem}
The goal of the adversarial problem is to maximize the infeasibility. In this problem, we consider the investment variables $\gamma_{ijt},\omega_{ijt},loi_{ijkt}$ and $\sigma_{it}$ as fixed since we want to evaluate the current investment solution obtained by solving the main problem at the current iteration. As mentioned in the previous section, we consider two types of infeasibilities: lack of generation capacity and line thermal rating violation. We consider them separately in consecutive problems. 

\paragraph{Generation infeasibility}
In this subproblem, we look for the random variable values that maximize the generation infeasibility. For that, we rewrite constraints \ref{active_balance} and \ref{reactive_balance} by including active and reactive power shedding respectively defined such as $P_{shed,it}\geq 0$ and $Q_{shed,it}\geq 0$.
\begin{flalign}
P_{Git}-\tilde{P}_{Cit}+P_{shed,it}=\sum_{\left(i,j\right)\ }{p_{ijt}} &
\label{active_balance_var}
\end{flalign}
\begin{flalign}
Q_{Git}-\tilde{Q}_{Cit}+Q_{shed,it}=\sum_{(i,j)\ }{q_{ijt}} &
\label{reactive_balance_var}
\end{flalign}
The adversarial subproblem corresponding to generation infeasibility is then written as follows:
\begin{flalign*}
\max_{\mathbf{\tilde{P}_{c}} \in \Omega^{P},\mathbf{\tilde{Q}_{c}} \in \Omega^{Q},P_{Git},Q_{Git}}\sum_{t=1}^T\sum_{i=1}^n A^{Gen}_{it}(P_{shed,it}+Q_{shed,it}) &\\
s.t. (12)-(13),(18)-(27)
\label{max_gen_infeas}
\end{flalign*}
$\mathbf{A^{Gen}}$ $\in \ \{0,1 \}^{n\times T}$ is simply a matrix that controls the indices $i$ and $t$ we want to include in the generation infeasibility maximization objective.

\paragraph{Line thermal rating infeasibility}
We now look for the random variable values that maximize the line thermal rating infeasibility. To this end, we remove constraint (23) for all indices and consider it in the objective. The corresponding adversarial subproblem is then:
\begin{flalign*}
\max_{\mathbf{\tilde{P}_{c}} \in \Omega^{P},\mathbf{\tilde{Q}_{c}} \in \Omega^{Q},P_{Git},Q_{Git}}\sum_{t=1}^T\sum_{ij} A^{Therm}_{ijt}(p_{ijt}^2+q_{ijt}^2-\gamma_{ijt}^2\overline{S}^2) &\\
s.t. (12)-(15),(18)-(22),(24)-(25)
\end{flalign*}
Similarly to the generation infeasibility adversarial problem, $\mathbf{A^{Therm}}$ $\in \ \{0,1 \}^{n\times n \times T}$ controls the indices $(i,j)$ and $t$ we want to include in the line thermal rating infeasibility maximization objective.
\subsection{Corrective problem}
In case the adversarial problem finds a problematic scenario, i.e. a scenario for which the objective of the adversarial problem is strictly positive, we now try to find corrective actions that can relieve the constraints violations previously maximized. The random variables are fixed (i.e. we fix the scenario) and we look for active/reactive generation setpoints such as to minimize constraint violations.
\paragraph{Corrective problem for generation infeasibility}
The problem is written as follows:
\begin{flalign*}
\min_{P_{Git},Q_{Git}}\sum_{t=1}^T\sum_{i=1}^n P_{shed,it}+Q_{shed,it} &\\
s.t. (12)-(13),(18)-(27)
\end{flalign*}
\paragraph{Corrective problem for line thermal rating infeasibility}
Eq. (23) is rewritten with a slack term $\delta_{ijt}\geq 0$ that represents a potential line thermal rating violation. It has to be noted that $\gamma_{ijt}$ is fixed in this problem, the introduction of $\delta_{ijt}$ in (23) thus doesn't remove its convexity. The problem is written as follows:
\begin{flalign*}
\min_{P_{Git},Q_{Git}}\sum_{t=1}^T\sum_{ij}\delta_{ijt}^2 \\
s.t. (12)-(15),(18)-(22),(24)-(25)&\\
p^2_{ijt}+q^2_{ijt}\le {\gamma }^2_{ijy}{\overline{S}}^2+\delta_{ijt}^2
\end{flalign*}
\subsection{Robust planning}
We present the whole algorithm for robust planning in the following flowchart (Fig \ref{flowchart}). For sake of clarity, we only describe the steps of the algorithms related to generation infeasibility. However, at every iteration of the algorithm, exactly equivalent steps are performed in parallel regarding line thermal rating infeasibility. Consequently, at every iteration, scenarios producing generation infeasibilities as well as scenarios causing line rating infeasibility are added to $\mathcal{S}$.
\tikzstyle{decision} = [diamond, draw,  
text width=6em, text badly centered, node distance=3cm, inner sep=0pt]
\tikzstyle{block} = [rectangle, draw, 
minimum width=9em, text centered, rounded corners, minimum height=4em]
\tikzstyle{block2} = [rectangle, draw, fill=yellow!20, 
text width=9em, text centered, rounded corners, minimum height=4em]
\tikzstyle{line} = [draw, -latex']
\tikzstyle{cloud} = [draw, ellipse,fill=red!20, node distance=3cm,
minimum height=4em]

\begin{figure}[!h]
\scalebox{0.65}{
  \centering
  \begin{tikzpicture}[node distance = 2cm, auto, scale=0.6]
    \node [block,yshift=-5cm] (node1) {\begin{tabular}{l} $\mathcal{S}\leftarrow \{\text{Deterministic scenario} \}$ \\ $\mathcal{PS}\leftarrow \emptyset$ \end{tabular}};
    \node [block, below of = node1] (node2) {\begin{tabular}{l} Unfix $\gamma,\omega,loi,\sigma$ \\ Solve main problem on $\mathcal{S}$\end{tabular}};
    \node [block, below of = node2] (node3) {\begin{tabular}{l} Fix $\gamma,\omega,loi,\sigma$\\$\mathcal{PS}\leftarrow \emptyset$ \end{tabular}};
    \node [block, below of = node3] (node4) {\begin{tabular}{l} $\forall (i^*,j^*), i^*\in \{1,...,n\}, t^* \in$ \{1,...,T\}\\ $A_{it}=1$ if $(i,t)=(i^*,t^*)$ and $0$ otherwise  \end{tabular}};
    \node [block, below of = node4] (node5) {\begin{tabular}{l} Solve generation adversarial problem\\Identify the set of violated constraints $\mathcal{VC}$\\$\mathcal{VC}=\{(i,t): Q_{shed,it}+P_{shed,it}> 0\}$ \end{tabular}};
    \node [block, below of = node5, yshift=-4.3cm] (node7) {\begin{tabular}{l} Fix $\tilde{P}_c$ and $\tilde{Q}_c$ and let $s^*=(\tilde{P}_c,\tilde{Q}_c)$\\Solve generation corrective problem \end{tabular}};
    \node [block, right of = node7, xshift=2.5cm,yshift=3.7cm] (node6) {\begin{tabular}{l}$A_{it}=1$ if $(i,t)\in \mathcal{VC}$\\ and $0$ otherwise \end{tabular}};
     \node [block, below of = node6] (node6bis) {\begin{tabular}{l}Solve generation \\adversarial problem \end{tabular}};
     \node [block,below of = node7,xshift=2cm] (node8) {$\mathcal{PS} \leftarrow \mathcal{PS}\cup \{ s^*\}$};
     \node [decision,below of = node8,xshift=-2cm] (node9) {All constraints considered?};
     \node [decision,below of = node9] (node10) {$\mathcal{PS} = \emptyset$?};
     \node [block,left of = node10,xshift=-2.5cm] (node12) {$\mathcal{S} \leftarrow \mathcal{S} \cup \mathcal{PS}$};
     \node [block,right of = node10,xshift=2cm] (node11) {END};
    \path [line] (node1) -- (node2);
    \path [line] (node2) -- (node3);
    \path [line] (node3)-- (node4);
    \path [line] (node4) -- (node5);
    \path [line] (node5) -- node{$\mathcal{VC}\supset \{(i^*,j^*)\}$}(node6);
    \path [line] (node6) -- (node6bis);
    \path [line] (node6bis) -- (node7);
    \path [line] (node7) -- node[xshift=0.5cm,yshift=-0.3cm]{$Objective>0$}(node8);
    \path [line] (node7) -- node[left]{$Objective=0$}(node9);
    \path [line] (node8) -- (node9);
    \path [line] (node12) |- (node2);
    \path [line] (node10) -- node{N} (node12);
    \path [line] (node10) -- node{Y} (node11);
    \path [line] (node9) -- node{Y} (node10);
    \path [line] (node9) --node{N}(-6.5,-40.5)--(-6.5,-18.5) |- (node4);
    \path [line] (node5) -- node{$\mathcal{VC}= \{(i^*,j^*)\}$}(node7);
    \path [line] (node5) --node{$\mathcal{VC}= \emptyset$}(-5.8,-25) -|(-5.8,-38)--(node9);
   
  \end{tikzpicture}
  }
  
  \caption{Flowchart of the robust planning algorithm for the case of generation infeasibility}
  \label{flowchart}
\end{figure}
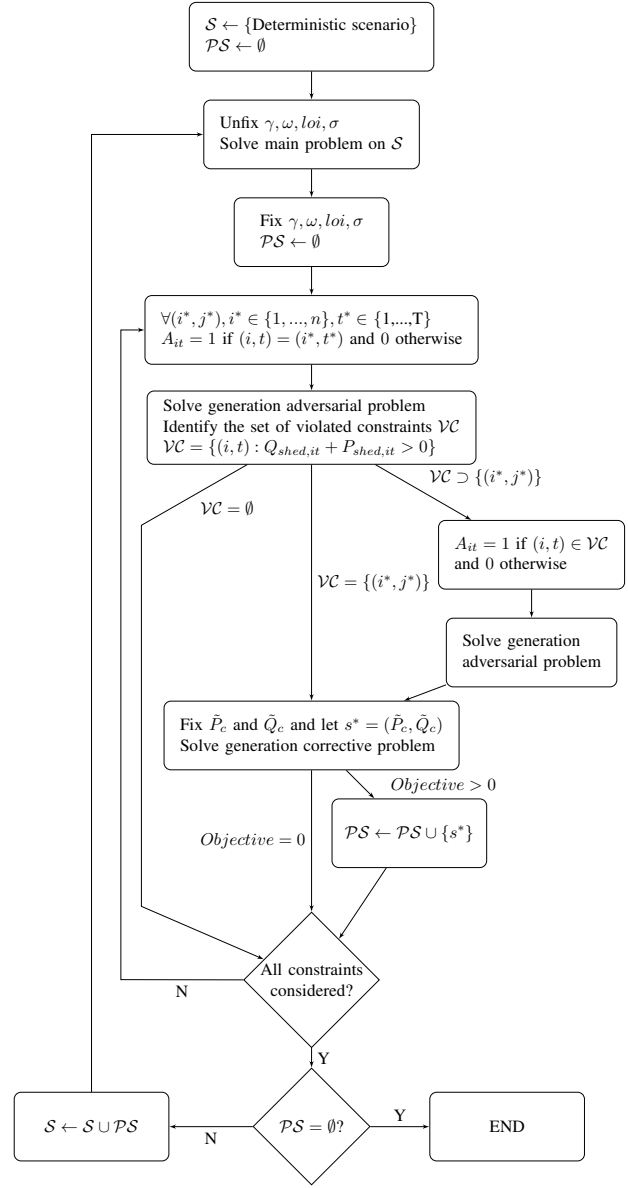
\subsection{Extension to general probabilistic modelling of uncertainty}
As mentioned in the beginning of this section, we only consider a rectangular uncertainty set with no assumption on the distribution of random variables. However, the proposed method can also be used with probabilistic modelling, i.e. when we consider that the joint distribution function of random variables is known (variables may be correlated in general). Indeed, let us consider the vector of random variables $\mathbf{\omega} \in \Omega$ and the joint density function $p(\mathbf{\omega})$. If we define the two vectors of parameters $\overline{\mathbf{\omega}}^L$ and $\overline{\mathbf{\omega}}^U$, the probability that $\mathbf{\overline{\omega}}^L \leq \mathbf{\omega}\leq\mathbf{\overline{\omega}}^U$ is then expressed as the following integral:
\begin{equation}
\mathbb{P}(\mathbf{\overline{\omega}}^L \leq \mathbf{\omega}\leq\mathbf{\overline{\omega}}^U)=\int_{\mathbf{\overline{\omega}}^L}^{\mathbf{\overline{\omega}}^U}p(\mathbf{\omega})d\omega
\end{equation}
As mentioned in \cite{margellos_road_2014}, this allows to formulate chance-constrained optimization as robust optimization. Indeed, the chance-constrained paradigm consists of finding the extremum of an objective function $f(x) $ while allowing constraints $h(x)\leq 0$ to be violated with a small probability $\epsilon$ as written hereunder. 
\begin{flalign}
&\sup_{x}f(x) \\
&s.t. \quad \mathbb{P}\big(h(x,\omega)\leq 0 \big)\geq 1-\epsilon \quad \forall \mathbf{\omega} \in \Omega
\label{chanceconstrained1}
\end{flalign}
We can reformulate this problem as a robust (deterministic) on a subspace of $\Omega$ such as the probability that random variables belong to this subspace is equal to $1-\epsilon$. This is written as follows. Note that the rectangular uncertainty interval defined by constraint (34) can be computed 'offline.
\begin{flalign}
\sup_{x}f(x)& \\
s.t.& \quad h(x,\omega)\leq 0 \quad \\
&\quad\mathbf{\overline{\omega}}^L \leq \mathbf{\omega}\leq\mathbf{\overline{\omega}}^U\\
&\quad\int_{\mathbf{\overline{\omega}}^L}^{\mathbf{\overline{\omega}}^U}p(\mathbf{\omega})d\omega = 1-\epsilon
\label{chanceconstrainedrobust}
\end{flalign}
\section{Results}
\label{Res}
The approach described in the previous sections is applied to a 20-node case described in \cite{carrano_electric_2006} on a 1-year planning horizon. Data for lines and loads consumptions for this case can be found in \cite{carrano_data_????-1}. Hourly consumption patterns are generated using real measurements used in \cite{navarro-espinosa_data_2015}. One representative day is considered for the whole year with 15 hourly consumption data. We consider a unique size for generators (2MW) and up to two lines placed in parallel between two nodes. The different models have been run on a 3.4Ghz Intel Core i7 processor with 8Go of memory. The models are written in AMPL and solved with CPLEX 12.7 using benders decomposition.\\ We compare the deterministic case where load consumptions $P_{cit}$ and $Q_{cit}$ are fixed and the uncertain case where they may vary between 50\% and 150\% of the deterministic value: $\tilde{P}_{cit}\in[0.5P_{cit};1.5P_{cit}]$ (idem for $\tilde{Q}_{cit}$). The results are shown in Table \ref{table1}.
\begin{table}[h!]
\caption{Comparison of planning solution for deterministic base case and robust case}
\label{table1}
\centering
\begin{tabular}{|r|c|c|}
\hline
&\small{Base case}&\small{Robust case}\\
\hline
\small{OPEX[M\$]}&\small{0.19}&\small{0.27}\\
\hline
\small{CAPEX[M\$]}&\small{3.30}&\small{4.02}\\
\hline
\small{Total cost[M\$]}&\small{3.49}&\small{4.29}\\
\hline
\small{Total amount of scenarios}&\small{/}&\small{274}\\
\hline
\small{Number of iterations}&\small{/}&\small{2}\\
\hline
\small{Computation time[s]}&\small{0.34}&\small{5079}\\
\hline
\end{tabular}
\end{table}
These results show that including load uncertainty in our planning problem increases the cost by almost 20\% for this test case. Indeed, the planning solution for the deterministic case only has 6MW of installed generation capacity while the robust solution has 8MW of installed generation capacity. It can also be observed that the computation time dramatically increases with the robust approach. As a matter of fact, the determination of problematic scenarios implies to solve several thousands of adversarial/corrective problems. Furthermore, at the second iteration, the main problem includes 274 times more operational constraints than the deterministic case which makes it a much bigger problem to solve than its deterministic counterpart.
\section{Conclusion}
\label{Concl}
In this paper, we presented a robust second-order cone formulation for the planning of autonomous microgrids under load uncertainty. An interval representation of uncertainty was used and it was shown that this approach could also be used to formulate a chance-constrained version of the planning problem. Preliminary results prove that inclusion of load uncertainty in the problem significantly increases the overall cost of the system which indicates the need to include uncertainty in planning. Further research will also include uncertainty related to RES-based generation and contingencies should be included as well to deliver more realistic planning solutions.

\bibliography{bibliotheque}{}
\bibliographystyle{ieeetr}

\end{document}